# Role of (periodic as well as aperiodic) tessellations in contemporary composition. The cases of *Tesselles sonores* and *Le Chapeau à douze cornes* by Marisa Acuña

Maria Luisa Acuña Fuentes*, Édouard Thomas
* Université Paris 8. UFR Arts, Philosophy and Aesthetics, **marisaacuna.f@gmail.com**

## Abstract

Tessellations represent an important area of geometry, where given shapes are duplicated and assembled to cover a surface without gaps or overlaps. The recent discovery of a family of aperiodic monotiles, which includes the famous "Hat", has shaken the field. Composers have already utilised the visual representation of plane tilings in their musical compositions. This is exemplified by Tom Johnson through his use of Vuza's canons, as well as other historical structures, such as the talea and color in isorhythmic motets. Constructions of the "Hat" from elementary geometric polygons provide a different perspective, for example through the sound transformation of microtonal intervals, as seen in Marisa Acuña's piece *Le Chapeau à douze cornes*. Furthermore, in the composer's piece *Tesselles sonores*, two types of periodic tilings are found: a Penrose tiling and an Alhambra tile (Nasrid Bird), which form the basis of some Euclidean plane tilings. The discovery of new Fujita and Niizeki tilings based on Penrose tilings also shows a similarity to the phase shifting technique used by Steve Reich in pieces like *Clapping Music*, creating a musical geometric pattern.

## 1. Introduction

Tessellations (tilings of a plane) constitute one of the major subjects of study in the field of mathematics. Their uses are highly diverse, spanning geometry, analysis, algebra, and beyond; their applications are ubiquitous in both art and the sciences. The main characteristic of these tessellations is the repetition of basic geometric shapes without leaving gaps between the shapes and without overlaps among the different shapes. The most well-known tessellations of the usual Euclidean space are based on regular, periodic repetition (through rotations, translations, or axial symmetries) of one or more core motifs. There are also aperiodic tilings, which notably lack any invariance under translation, such as the famous Penrose tilings. These are constructed using a set of elementary polygonal shapes with appropriate matching rules or decorations to prevent them from tiling periodically. Between 1974 and 1976, Roger Penrose created sets containing six, then four, and finally two polygons (the *kite* and the *arrow*), enabling the generation of aperiodic tilings of the Euclidean plane (Gardner, 1977; Penrose, 1979). From that point on, a fundamental mathematical question remained open: was it possible to obtain an aperiodic tiling of the plane from a single hypothetical (connected, polygonal) shape, referred to as an *aperiodic monotile*?
This question received a complete (positive) answer in 2023 with the discovery of the "Hat" by David Smith, followed by the "Turtle," then the "Spectre," and their remarkable mathematical properties (Smith *et al.*, 2024*a*, 2024*b*); see Section 2. The "Hat" thus became the first known (compact, connected, polygonal) aperiodic monotile. This polygon can be constructed by anyone in a completely elementary manner (using sixteen isometric right triangles, or eight kite-shaped polygons, or three regular hexagons arranged to tile the plane); see Section 3. This immediate appropriation of a new mathematical object has contributed to a global enthusiasm, engaging artists, mathematicians, scientists, and the general public, as evidenced by the million views of Jade Tan-Holmes's video on her online channel Up and Atom dedicated to the "Hat" (Tan-Holmes, 2023) and by the plethora of science enthusiasts who created, recorded and shared mathematical contents of the highest quality (see Haran 2023*a*, 2023*b*, 2023*c*; Hearn, 2023; MacDonald, 2023*a*, 2023*b*; Trappenberg, 2023 among many others).
To illustrate this enthusiasm, we will present an example of the use of the properties of aperiodic tilings in the construction of a contemporary musical work.

## 2. Tilings and aperiodicity

The discovery of the "Hat" in 2023 by David Smith marks a historic mathematical event in the study of aperiodic tilings of the Euclidean plane. Smith uncovered a basic polygonal geometric figure (named the "Hat", presumably due to its shape) capable of generating a tiling of the plane that possesses neither translational nor rotational symmetry, making the "Hat" the first known example of an aperiodic monotile. An infinite number of other polygonal aperiodic monotiles have been obtained simply by deforming the "Hat", all belonging to the same family (Smith *et al.*, 2024*a*).
However, some instances of the "Hat" must be reflected using axial symmetry to create an aperiodic tiling of the plane; a semantic controversy arose regarding whether the resulting aperiodic tiling was obtained using a single shape (the "Hat") or two distinct shapes (the "Hat" and its mirror image). By exploiting the properties of the "Spectre" (another aperiodic monotile from the same family as the "Hat"), this controversy was resolved: thanks to the "Spectre", one can achieve a true unique aperiodic monotile without needing to use a mirror image of that aperiodic monotile. For the musical development discussed in this article, we can, however, operate under the assumption that the "Hat" is indeed an aperiodic monotile: the notion of a mirror image is not relevant in our application.

The "Hat" is a polygon with thirteen sides (see Figure 7). It consists of eight isometric kites; it is therefore a poly-kite (Smith *et al.*, 2024*a*). The musical analysis of Marisa Acuña's piece *Le Chapeau à douze cornes* will highlight various geometric transformations involved in constructing the "Hat".
Since then, several authors have explored in greater detail the mathematical and structural properties of the "Hat" and the family of aperiodic monotiles that derive from it. For example, Socolar (2023) examines the quasicrystalline structure that emerges from these tilings, revealing its complexity and its relationship to the characteristic patterns of quasicrystals. In a complementary spirit, Baake (2025) studies the dynamics and topology of the family of tilings generated by the "Hat", highlighting the local symmetries and global characteristics of these structures. Meanwhile, Imura (2025) proposes a new family of non-periodic tilings that, using elementary geometric tools, present a type of structural regularity that allows the "Hat" to be placed in the broader context of recent research on aperiodicity. Fujita and Niizeki's (2025) work on a new variant of the Penrose rhombus tiling brings new insights to the study of non-periodic tilings. This article presents a geometric adaptation that preserves the characteristic aperiodicity of the Penrose pattern while introducing elements of formal simplicity that may facilitate its analysis or application to other fields, such as music.

On a purely mathematical level, research continues to seek other examples of aperiodic monotiles that do not belong to the family uncovered by David Smith and his co-authors. In Jedrzejewski's (2009) research, the concept of aperiodic tiling is not limited to geometric planes and extends its scope to the set of integers (or to a group G, finite or infinite), thereby making direct connections with the structure of sequences and rhythmic patterns over time. This algebraic approach, currently flourishing due to the spectacular advances of Tao and Greenfeld (2023), provides a particularly relevant framework for the analysis and creation of musical systems based on aperiodicity.

In this article, we will attempt to connect these theories to musical works where they may be present, consciously or unconsciously, in the composer's work.

## 2.1. Music and tilings: previous work

In recent decades, many composers and theorists have explored the possibilities offered by tilings as structural tools in musical composition. Whether through a rhythmic, spatial, or purely mathematical approach, this research has laid the groundwork for a rigorous and creative interdisciplinary field between music and geometry. For an authoritative introduction to tessellations that is suitable for both students, professional mathematicians, artists, and craftsmen, see Graunbaum and Shepherd (2016). For a short and pleasant introduction to how the "Hat" was discovered by David Smith and the collaboration that followed, see Kaplan (2024), where you will also find a very accessible overview of some mathematical challenges still waiting to be addressed in tiling theory.

One of the pioneers to have explicitly addressed this relationship between music and geometry is American composer Tom Johnson (1939−2024), who describes in his essay *Tiling in My Music* (Tom Johnson, 2003) how he utilised tiling patterns —particularly through rhythmic and spatial sequences— as a compositional principle in several of his works. In his approach, repetition, the precise arrangement of motifs, and local symmetries reflect the logic of geometric tilings. Johnson references the studies of Romanian mathematician Dan Tudor Vuza (Vuza, 1991, 1992*a*, 1992*b*, 1993) (see also Andreatta, 2011), who investigates how a rhythmic subset R of the finite group $\mathbb{Z}/n\mathbb{Z}$ can be "completed" by another subset S of $\mathbb{Z}/n\mathbb{Z}$ in order to "fill the tiling" (cover the entire rhythmic structure without overlaps). Specifically, in the piece *Tilework for Clarinet* (Tom Johnson, 2011), the composer utilises the group $\mathbb{Z}/15\mathbb{Z}$ and two elements: the original rhythm R = {0, 2, 5} (which means that rhythmic strikes occur at times 0, 2, and 5), and a rhythm associated with rhythm R, namely R′ = {0, 3, 5}.

- The clarinet enters several times, repeating either the original rhythm {0, 2, 5} or what Tom Johnson calls its *retrograde*, namely the rhythm R′ = {0, 3, 5};
- The entries are staggered in time, as in a canon;
- Each entry fills in the "gaps" left by the previous entries.

An example of layout (see also Figure 1):

| **Voices** | **Rhythm used** | **Begins at time** |
|---|---|---|
| Voice 1 | (0, 2, 5) | 1 |
| Voice 2 | (0, 2, 5) | 2 |
| Voice 3 | (0, 3, 5) | 5 |
| Voice 4 | (0, 3, 5) | 12 |
| Voice 5 | (0, 2, 5) | 9 |
| Voice 6 | (0, 3, 5) | 13 |

Each voice is a translation of one of the two motifs R or R′. If $t \in \mathbb{Z}/15\mathbb{Z}$ represents the entry time of a voice, then:

- a voice with R will play R + $t$, namely $\{(r + t) \bmod 15, \text{ with } r \in \mathrm{R}\}$;
- a voice with R′ will play R′ + $t$, namely $\{(r + t) \bmod 15, \text{ with } r \in \mathrm{R}'\}$.

Entry time of the voices (see also Figure 1):

| **Voice** | **Motif used** | **Entry** |
|---|---|---|
| 1 | R | 1 |
| 2 | R | 2 |
| 3 | R’ | 5 |
| 4 | R’ | 12 |
| 5 | R | 9 |
| 6 | R’ | 13 |

We define the set $T$ of all the times utilised: $T = \bigcup_{i=1}^{6} (\mathrm{R}_i + t_i) \bmod 15$, where $\mathrm{R}_1 = \mathrm{R}_2 = \mathrm{R}_5 = \mathrm{R}$ and $\mathrm{R}_3 = \mathrm{R}_4 = \mathrm{R}_6 = \mathrm{R}'$ and $t_i$ is the entry time of the voice $i$. We then get:

$T = \{0, 1, 2, 3, 4, 5, 6, 7, 8, 9, 10, 11, 12, 13, 14\}$, which can be identified to $\mathbb{Z}/15\mathbb{Z}$.

Formally, we have R = {0, 2, 5} and R′ = {0, 3, 5}; let B = {1, 2, 9} and B′ = {5, 12, 13} be their respective entry times (see also Figure 1). We then have:

$$\bigcup_{b \in B} (\mathrm{R} + b) \cup \bigcup_{b' \in B'} (\mathrm{R}' + b') \cong \mathbb{Z}/15\mathbb{Z}.$$

This composition retains a formal resemblance, in its construction, to the isorhythmic motets used in the Middle Ages.

### 2.2. Presentation of the formal mathematical model

Consider the following elements defined over $\mathbb{Z}$, the set of relative integers (or possibly over a finite group $\mathbb{Z}/n\mathbb{Z}$ if one wishes to modulate the time):

#### A-Talea (rhythmic structure)

Let $\boldsymbol{T} = (t_1, t_2 \ldots t_m)$ be a vector of rhythmic durations (which are all whole numbers), with $t_i > 0$ for all integer $i$ between 1 and $m$. This sequence repeats itself in a cyclic way.

Example:

$\boldsymbol{T} = (2, 1, 1)$, namely quarter note, eighth note, eighth note.

#### B-Color (melodic structure)

Let $\boldsymbol{C} = (c_1, c_2 \ldots c_n)$ be a sequence of pitches (or intervals), represented by relative integers (for instance: MIDI pitches, or distances from a central C). This sequence repeats itself also in a cyclic way.

Example:

$\boldsymbol{C}$ = (C, E, F, G), namely four distinct pitches. With MIDI values, we can write:
$\boldsymbol{C} = (60, 62, 64, 65)$.

- Talea: $\boldsymbol{T} = (2, 1, 1)$, $\boldsymbol{T} = (2, 1, 1)$, $\boldsymbol{T} = (2, 1, 1)$, with $m = 3$;
- Color: $\boldsymbol{C} = (60, 62, 64, 65)$, $\boldsymbol{C} = (60, 62, 64, 65)$, $\boldsymbol{C} = (60, 62, 64, 65)$, with $n = 4$;
- N = LCM(3, 4) = 12.

We next introduce $M = \left( C_{(i \bmod n)+1} \quad T_{(i \bmod m)+1} \right)_{0 \le i \le N-1}$ in order to store in a same matrix all the useful data. We easily determine every coefficient of $M$:

For $i = 0$: $C_{(0 \bmod 4)+1} = C_1 = 60$; $T_{(0 \bmod 3)+1} = T_1 = 2$;

For $i = 1$: $C_{(1 \bmod 4)+1} = C_2 = 62$; $T_{(1 \bmod 3)+1} = T_2 = 1$;

…

For $i = N - 1 = 11$: $C_{(11 \bmod 4)+1} = C_4 = 65$; $T_{(11 \bmod 3)+1} = T_3 = 1$.

Results:

| Step | Pitch (color) | Duration (talea) |
|---|---|---|
| 0 | 60 | 2 |
| 1 | 62 | 1 |
| 2 | 64 | 1 |
| 3 | 65 | 2 |
| 4 | 60 | 1 |
| 5 | 62 | 1 |
| 6 | 64 | 2 |
| 7 | 65 | 1 |
| 8 | 60 | 1 |
| 9 | 62 | 2 |
| 10 | 64 | 1 |
| 11 | 65 | 1 |

A graphic representation (see Figure 1) allows one to compare, visually, the movements. In the case of the Vuza's canon with $\mathbb{Z}/15\mathbb{Z}$:

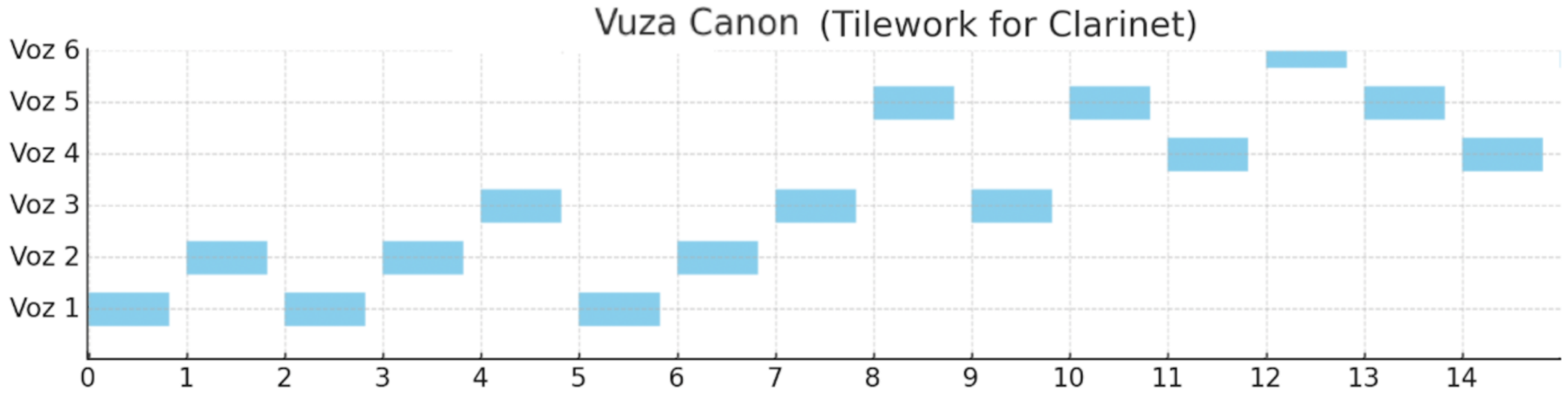

**Figure 1**. Representation of *Tilework for Clarinet* (Tom Johnson, 2011).

In the case of the isorhythmic motet with $\mathbb{Z}/12\mathbb{Z}$:

Isorhythmic Motet Diagram (Mod 12): Talea and Color

| T1 | T2 | T3 | T1 | T2 | T3 | T1 | T2 | T3 | T1 | T2 | T3 |
|---|---|---|---|---|---|---|---|---|---|---|---|
| A | B | C | D | A | B | C | D | A | B | C | D |

**Figure 2**. Diagram talea and color for the modulo 12.

As we have illustrated, the construction of the motet can bear some similarities to Vuza's canon in the sense that both seek to create a shifted sound and rhythmic model in the order of entry of the voices. In Johnson's piece, two rhythmic formulas create a (rhythmic) tiling through the staggering of entries across different voices. In the case of the isorhythmic motet, a shift is applied to the same voice between the color (notes) and the talea (rhythm). The isorhythmic motet does not aim to cover the entire rhythmic space. However −and as we will see in the following sections of this article− the shifting of a single rhythmic voice can give rise to another type of tiling.

## 3. Analysis of two Musical Compositions. *Le Chapeau à douze cornes* and *Tesselles sonores*.

### 3.1. *Le Chapeau à douze cornes*

The piece *Le Chapeau à douze cornes* ("The Twelve-Horned Hat") was composed for flute, violin, viola, and cello. There is a main difference between this musical composition and those based on a combinatorial use of elements or canon to describe a tessellation: on *Le Chapeau à douze cornes* (Marisa Acuña, 2025), the composer aimed to depict the geometric transformation that takes one from an elementary geometric figure (a regular hexagon) to the "Hat", a newly discovered aperiodic monotile. The intention is not to cover a space (rhythmic or otherwise) with a tessellation or to repeat a basic motif, but rather to illustrate the evolution of the transformation from a regular hexagon (which allows for a periodic tiling of the plane) into an aperiodic monotile (the "Hat").

Two ideas are explored on *Le Chapeau à douze cornes*: on one hand, an aesthetic idea inspired by the music of Manuel de Falla and flamenco; on the other hand, a geometric idea (which we will focus on), which forms the fundamental axis of the material used for the development of the piece.

The piece is divided into five movements. Each describes the geometric transformation of a reference regular hexagon H into the "Hat", accompanied by a transformation of the sound and rhythmic material, which we will now detail.

Movement I, "*Hexagone*" ("Hexagon"), begins with a six-note scale (notes selected according to the piece's aesthetic). Each note corresponds to one side of a regular hexagon H.

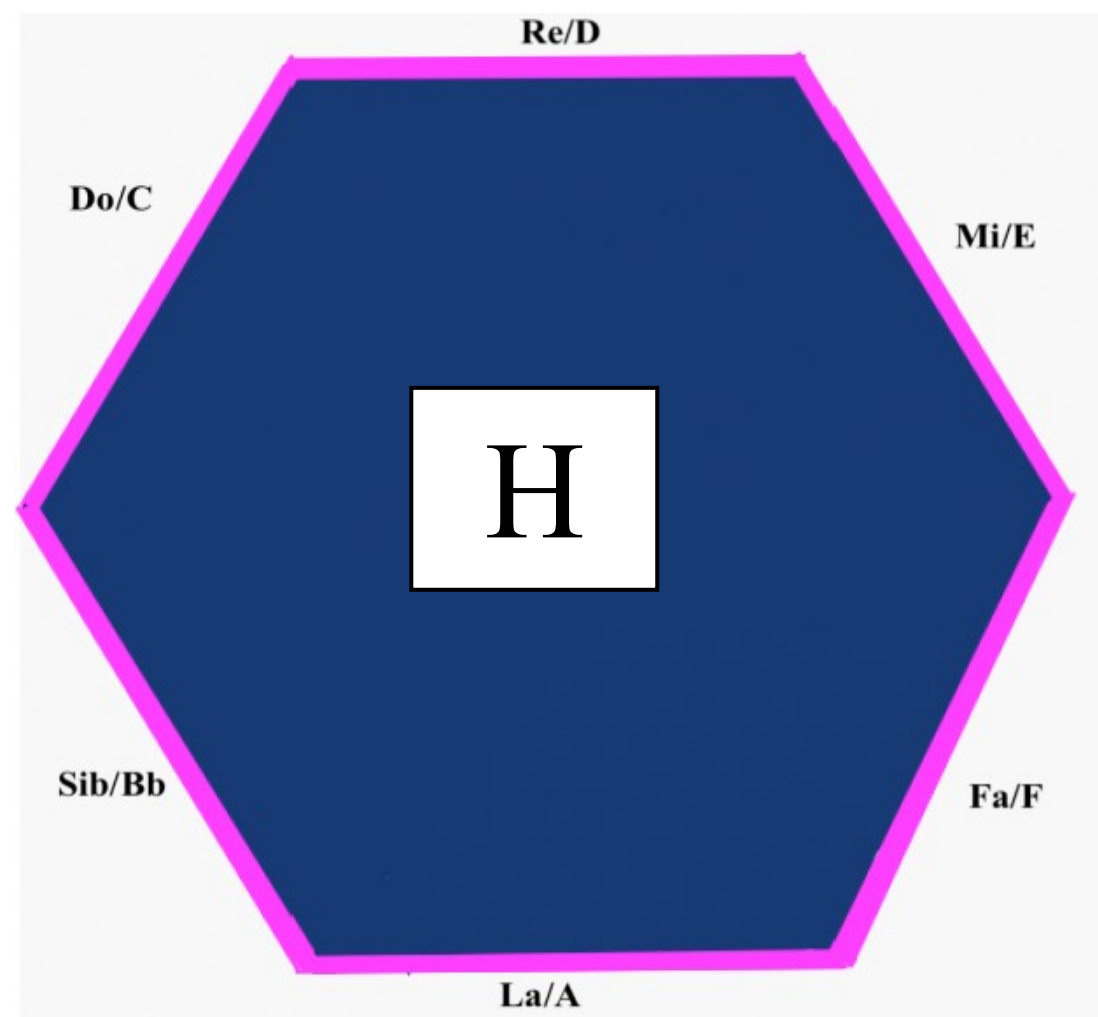

**Figure 3**. Notes on the regular hexagon H.

The transition from this movement to Movement II will be linked to the sounds and the addition of two six-note scales in order to depict the three hexagons from which the "Hat" will be constructed (see Figure 4). Thus, the initial idea of Movement II, "*Trois hexagones*" ("Three Hexagons"), is the presentation of the three isometric regular hexagons used, $H_1$, $H_2$ and $H_3$. Each hexagon, a copy of the reference regular hexagon H (with $H_1$ = H), uses a different scale. The violin will use the same scale as in the first movement (1), but transposed. The viola and flute will share the same scale (2). The cello will use a third, distinct scale (3).

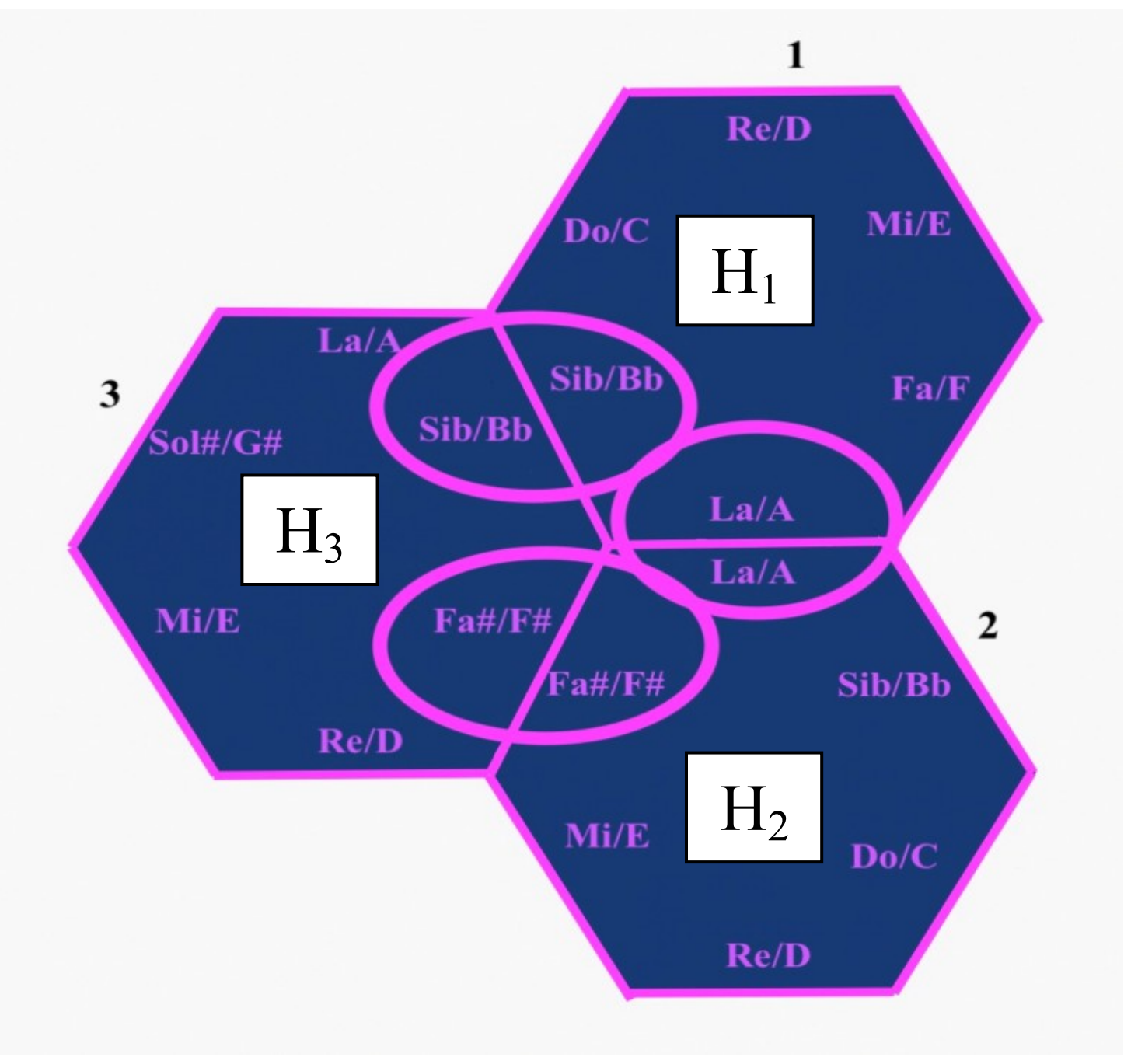

**Figure 4**. Three copies, $H_1$, $H_2$ and $H_3$, of the regular hexagon H (not to scale).

Three notes (circled in Figure 4) are associated with the three common sides of two hexagons: B♭, A, and F♯; they will be used throughout the rest of the movement. Let us illustrate these first two movements represented in *pitch class* values to facilitate the understanding of this transformation.

**Movement I:** “**Hexagon**”

Notes: A–B♭–C–D–E–F (see Figure 3).

Translation into pitch classes[1] in $\mathbb{Z}/12\mathbb{Z}$ (C = 0):

| Note | Pitch class (mod 12) |
| --- | --- |
| A | 9 |
| B♭ | 10 |
| C | 0 |
| D | 2 |
| E | 4 |
| F | 5 |

→ **Hexagon I ($H_1$)**:

Same as H (see Figure 4). Pitch classes values: {0, 2, 4, 5, 9, 10}.

**Movement II:** “**Three Hexagons**”

**Hexagon I ($H_1$)**

Pitch classes values: same as Movement I, {0, 2, 4, 5, 9, 10}.

**Hexagon II ($H_2$)**

| Note | Pitch class (mod 12) |
| --- | --- |
| C | 0 |
| D | 2 |
| E | 4 |
| F♯ | 6 |
| A | 9 |
| B♭ | 10 |

Pitch classes values (see Figure 4): {0, 2, 4, 6, 9, 10}.

[1] *Pitch class*: notation system in which each note of the chromatic scale is encoded by a number from 0 to 11. Thus C = 0, C# = 1…

**Hexagon III ($H_3$)**

Notes (see Figure 4): E–F♯–G♯–A–B♭–D.

| **Note** | **Pitch Class** |
|---|---|
| E | 4 |
| F♯ | 6 |
| G♯ | 8 |
| A | 9 |
| B♭ | 10 |
| D | 2 |

Pitch classes values: {2, 4, 6, 8, 9, 10}.

The rhythm itself is based on the idea of a motif derived from the number 6, beginning with percussion on the instrument's wood, characteristic of the rhythm of *bulería*, which is organized as follows:

**Accent vector (1 = accent, 0 = no accent):**

[1, 0, 0, 1, 0, 0, 1, 0, 1, 0, 1, 0].

| **Beat** | **1** | **2** | **3** | **4** | **5** | **6** | **7** | **8** | **9** | **10** | **11** | **12** |
|---|---|---|---|---|---|---|---|---|---|---|---|---|
| Stroke | ● | | | ● | | | ● | | ● | | ● | |

**Typical grouping: 6 (3 + 3) and 6 beats:**

(1 2 3) (4 5 6) (7 8 9 10 11 12).

Movement III, "*Cerf-volant*" ("Kite"), introduces a large, irregular, symmetric, non-convex hexagon in the shape of a kite-like polygon, formed by five of the eight isometric convex small kites located at the center of the three hexagons (see Figure 5). These eight small kites are convex quadrilaterals that naturally make up the "Hat"; they consist of two isometric right-angled triangles assembled together. Here, each side of the large kite-like shape is either half a side of the hexagon (for two sides), an apothem of the hexagon (for two other sides), or two apothems from two different hexagons (for the remaining two sides). In total, there are two half-sides (from two hexagons) and six apothems (two from each hexagon). Finally, it should be noted that none of the small kites and the large kite-like shape discussed in this article are related to the famous "kites" introduced along the "arrows" in the 1970s by Penrose in his aperiodic tilings (see Introduction).

The construction of the large kite (in green on the sketch below) starts from one side of a hexagon, which transforms into a half-side: if the initial note moves to a higher pitch, it will be raised by a semitone; if it moves to a lower pitch, it will be lowered by a semitone.

If we fix the length of a half-side of a hexagon to 1, some elementary geometry allow us to determine the lengths of the sides formed by the apothems. The actual value of an apothem is √3, which is an irrational number. To make this value feasible for musicians, the composer approximated the apothem to three-quarters of the side of the hexagon (the full side of the hexagon equals two half-sides). In fact, to within $10^{-2}$, √3 is approximately 1.73, while ¾ × 2 is 1.5. Therefore, if the side of the hexagon corresponds to a whole tone, the apothem will be considered as three-quarters of a tone. These ¾ tones will be added or subtracted to the starting note depending on the direction the large kite takes toward the next note. To sum up, the encoding of magnitudes is as follows: 2 for a whole tone; 1 for a semitone; 1.5 for three quarters of a tone. Therefore, the use of microtonal intervals (smaller than a semitone, so including quarter-tone long ones) in the piece is not merely the result of an arithmetic subdivision of tones. Rather, it emerges directly from the geometric transformations involved in the transition from the periodic regular hexagon to the aperiodic monotile. The appearance of irrational geometric proportions in the construction of the "Hat" naturally suggested the use of microtonal pitch organization. In this way, the three hexagons are transformed, on a sonic level, into a large kite-like shape (in green in Figure 5).

The central section, the most rhythmic section of this movement, alternates between 3/4 and 2/4 time signatures in a sequence of 3 measures (three beats) + 2 measures (two beats) + 3 measures (three beats), where the 2/4 measures represent the two half-sides, and the 3/4 measures symbolize the appearance of the apothems.

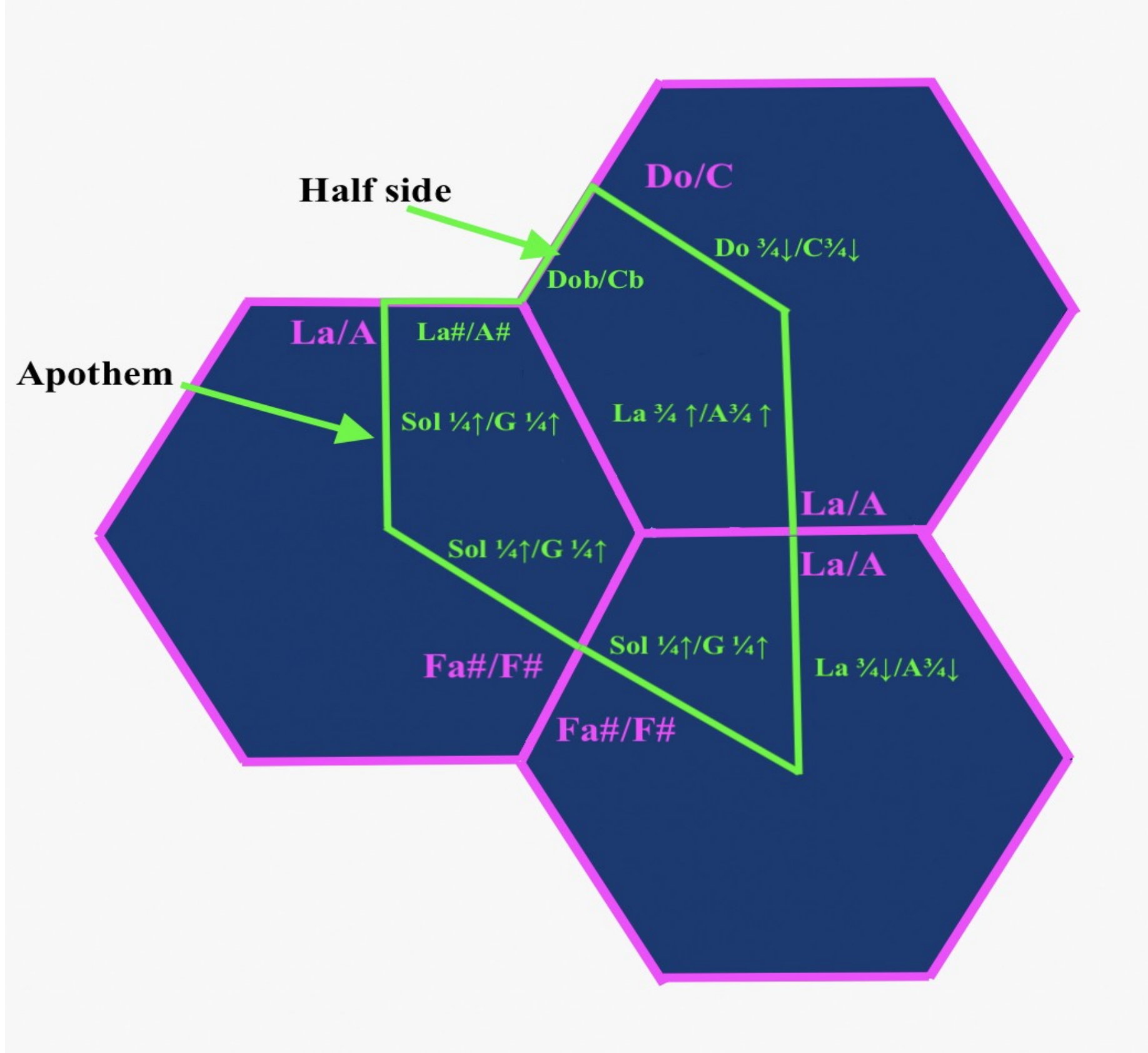

**Figure 5**. In green, the large kite on three hexagons.

## Geometric representation of the contour

We can represent the contour of the large kite-like shape as a sequence of six vectors in the plane: $v_1$, $v_2$, $v_3$, $v_4$, $v_5$, $v_6$.

Each vector has a *direction* (angle with respect to a fixed axis) and a *magnitude* (length). For the magnitude, the following encoding was used: 2 for the whole tone, therefore 1 for the semitone and 1.5 for three quarters of a tone (intervals smaller than a semitone, so including quartertone long intervals, are called microtonal intervals in music theory).

| Side of the *large kite* | Type | Length (approx.) | Musical interval |
|---|---|---|---|
| $v_1$ | Apothem | ≈ 1.5 | ±¾ tone |
| $v_2$ | Apothem | ≈ 1.5 | ±¾ tone |
| $v_3$ | Apothem | ≈ 1.5 | ±¾ tone |
| $v_4$ | Half-side of hex. | 1 | ±½ tone |
| $v_5$ | Half-side of hex. | 1 | ±½ tone |
| $v_6$ | Two apothems | ≈ 3.0 | ±1½ tone |

| Step | Base note | Transformation (semitones) | Result (mod 12) |
|---|---|---|---|
| 1 | C (0) | −1 | 11 (C♭) |
| 2 | A (9) | +1 | 10 (A#) |
| 3 | C (0) | −1.5 | 10.5 (C − 3/4) |
| 4 | A (9) | +1.5 | 10.5 (A +3/4) |
| 5 | A (9) | −1.5 | 7.5 (A − 3/4)(G+1/4) |
| 6 | F♯ (6) | +1.5 | 7.5 (G +1/4) |
| 7 | F♯ (6) | +1.5 | 7.5 (G +1/4) |
| 8 | F♯ (6) | +1.5 | 7.5 (G +1/4) |

Movement IV, "*Triangles isometriques*" ("Isometric Triangles"), is made up of the eight small kites that form the "Hat"; these are quadrilaterals but the composer preferred to call them "isometric triangles" because each of these eight quadrilaterals is made up of two true isometric right-angled triangles. In this movement, each small kite consists of two half-sides and two apothems. The development of the sonic material follows the same principle as in Movement III. The material is organized as follows: four small kites belong to the first hexagon, and the remaining four are distributed in pairs across the other two hexagons. Thus, there is a sequence of four measures using the first material for the four first small kites; then two measures with the material of one hexagon; then two measures with the material of the last hexagon. Once the sequence (4 + 2 + 2) is completed, it repeats. The rhythm alternates between measures of six beats and three beats.

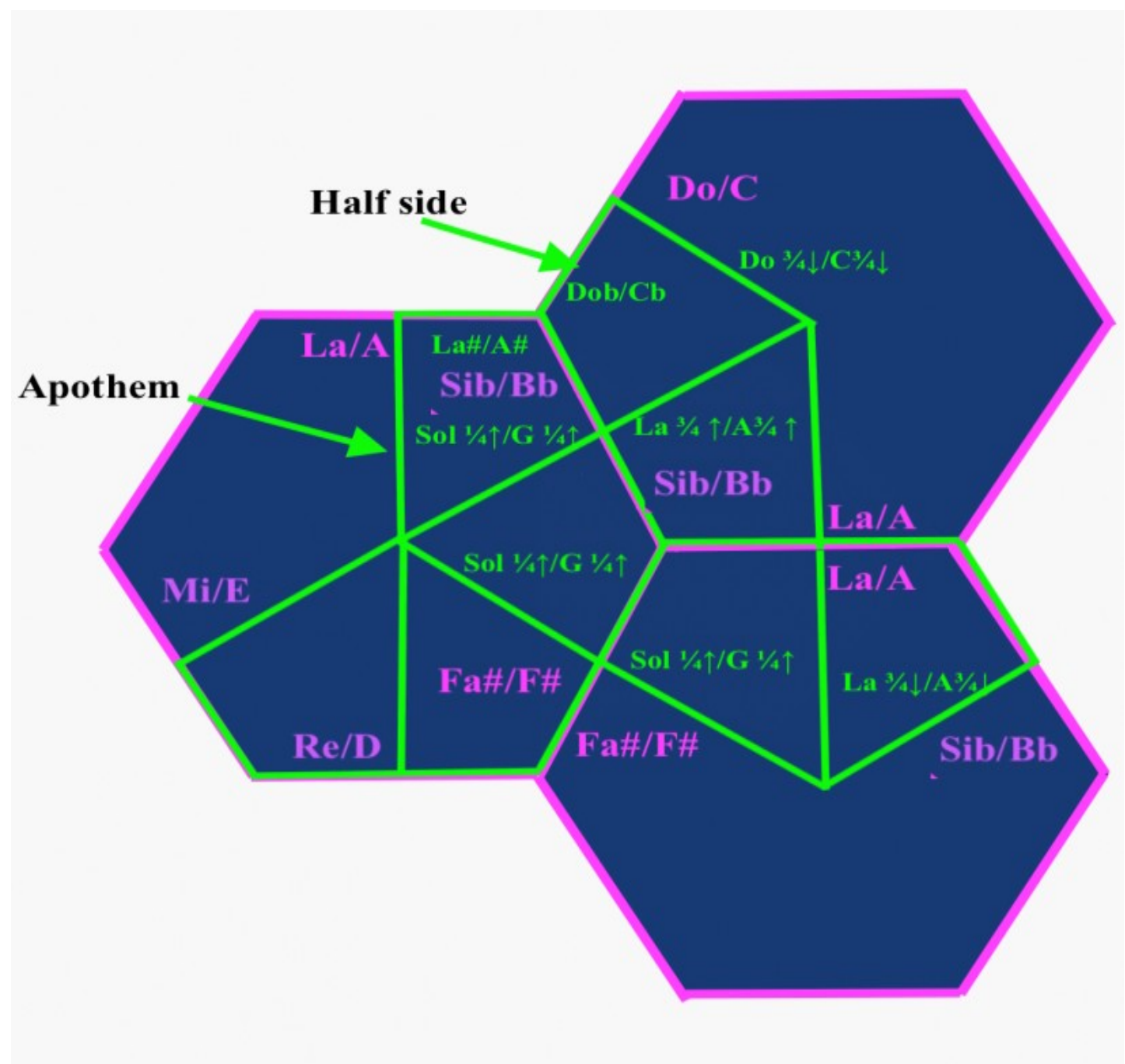

**Figure 6**. In green, the eight small kites on three hexagons that the composer chose to name "isometric triangles" in *Le Chapeau à douze cornes* (Marisa Acuña, 2025).

Each row corresponds to one of the eight small kites that make up the "Hat". The transformation vector of each small kite can be written as $v_i = [\pm1, \pm1, \pm1.5, \pm1.5]$, where:

- ±1 corresponds to a half-side, therefore to 1 semitone**;**
- ±1.5 corresponds to an apothem, therefore to 1.5 semitone (¾ of a whole tone).

The rhythm follows an alternating pattern of 6/8 and 3/8 measures.

The form follows a 4 + 2 + 2structure, with respect to the position of the eight small kites:

- Four small kites on the left of Figure 6 belong to hexagon $H_3$,
- Two small kites on the right upper corner of Figure 6 are located in $H_1$,
- Two small kites on the right lower corner of Figure 6 are part of $H_2$.

| Small kite (*i*) | Hexagon | Time signature (beats) | Transformations (in semitones) |
|---|---|---|---|
| 1 | $H_3$ | 6 | [−1, +1, −1.5, +1.5] |
| 2 | $H_3$ | 3 | [+1, −1, +1.5, −1.5] |
| 3 | $H_3$ | 6 | [−1, −1, +1.5, +1.5] |
| 4 | $H_3$ | 3 | [+1, +1, −1.5, −1.5] |
| 5 | $H_1$ | 6 | [−1, +1, +1.5, +1.5] |
| 6 | $H_1$ | 3 | [+1, −1, −1.5, −1.5] |
| 7 | $H_2$ | 6 | [−1, −1, −1.5, +1.5] |
| 8 | $H_2$ | 3 | [+1, +1, +1.5, −1.5] |

For Movement V, “*Le Chapeau*” (“The Hat”), the construction of the sonic material is the same as in Movements III and IV. Here, the notes will be determined by each side, half-side, or apothem of the three hexagons. For this final movement, the main rhythmic signature is 6.

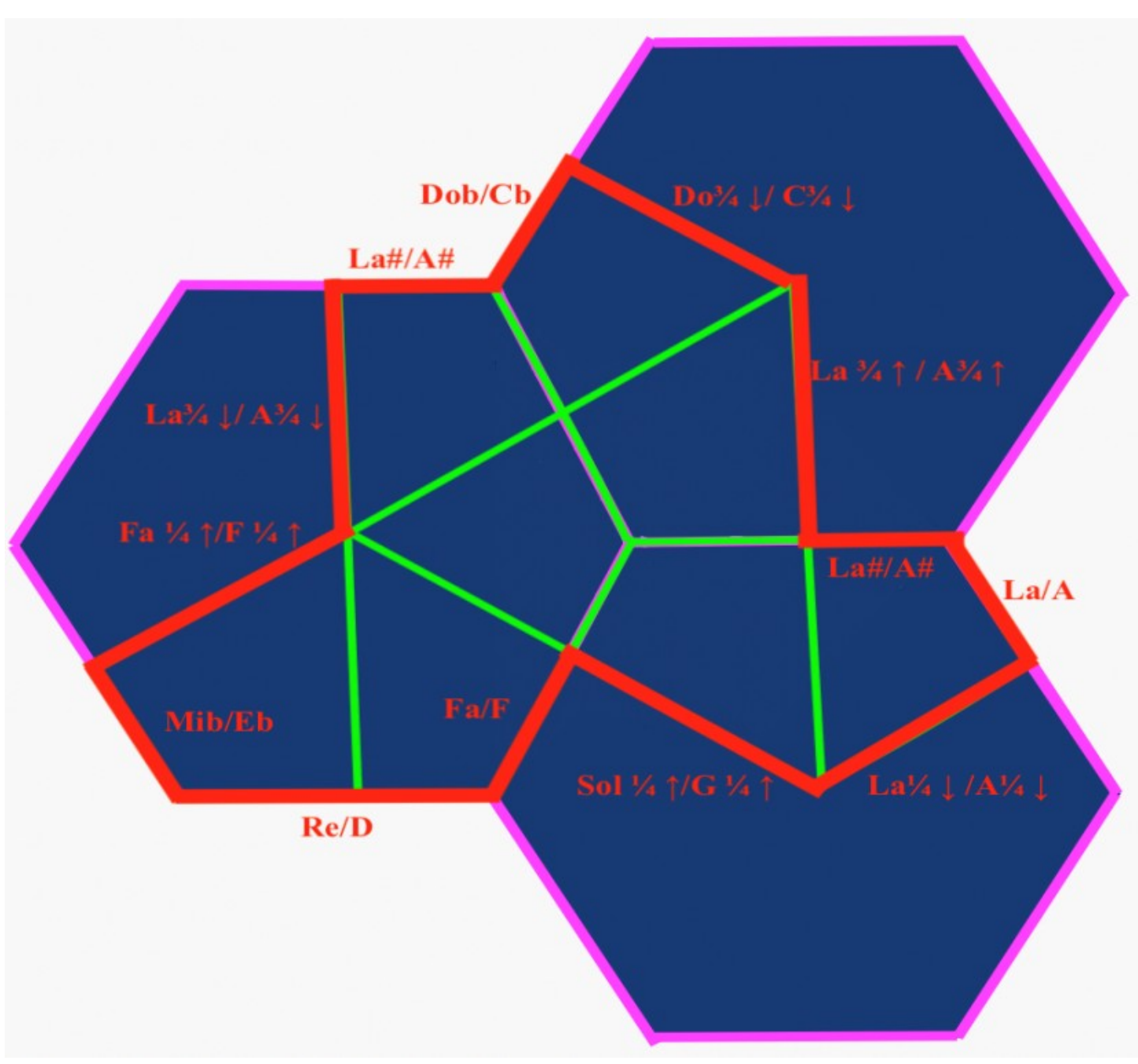

**Figure 7**. In red, the “Hat” is a 13-sided polygon defined by one full side, six half-sides and six apothems from three regular hexagons.

| Nº. | Segment type | Base note | Transformation (semitones) | Result (mod 12) |
|---|---|---|---|---|
| 1 | Half-side | C (0) | −1 | 11 |
| 2 | Half-side | A (9) | +1 | 10 |
| 3 | Apothem | A (9) | +1.5 | 10.5 |
| 4 | Apothem | E (4) | +1.5 | 5.5 |
| 5 | Half-side | E♭ (3) | −1 | 2 |
| 6 | Full side | D (2) | 0 | 2 |
| 7 | Half-side | F♯ (6) | −1 | 5 |
| 8 | Apothem | F♯ (6) | +1.5 | 7.5 |
| 9 | Apothem | B♭ (10) | −1.5 | 8.5 |
| 10 | Half-side | B♭ (10) | −1 | 9 |
| 11 | Half-side | A (9) | +1 | 10 |
| 12 | Apothem | A (9) | +1.5 | 10.5 |
| 13 | Apothem | C (0) | −1.5 | 10.5 |

### 3.2. *Tesselles sonores*

The musical piece *Tesselles sonores* (Marisa Acuña, 2022), written before the discovery of the first aperiodic monotile, uses models closer to the idea developed by the canon of Vuza in Tom Johnson's works (see Section 2). *Tesselles sonores* ("Sound Tessellation") is a three-movement piece written for mandolin, harp, and guitar.

Two types of tessellations were used to construct this piece. For the first movement, one of Penrose's tessellations was used. It consists of two different types of rhombuses (fat rhombuses, in blue, and thin rhombuses, in green) that can be assembled in several different ways. Here are the four combinations chosen by the composer for the piece:

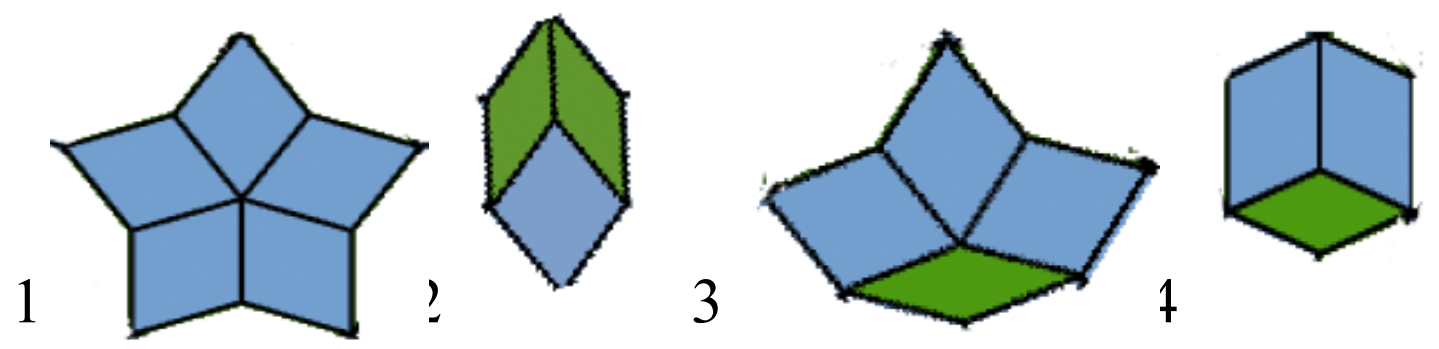

**Figure 8**. The four combinations (not to scale)
of the fat and thin rhombuses from Penrose's tessellation.

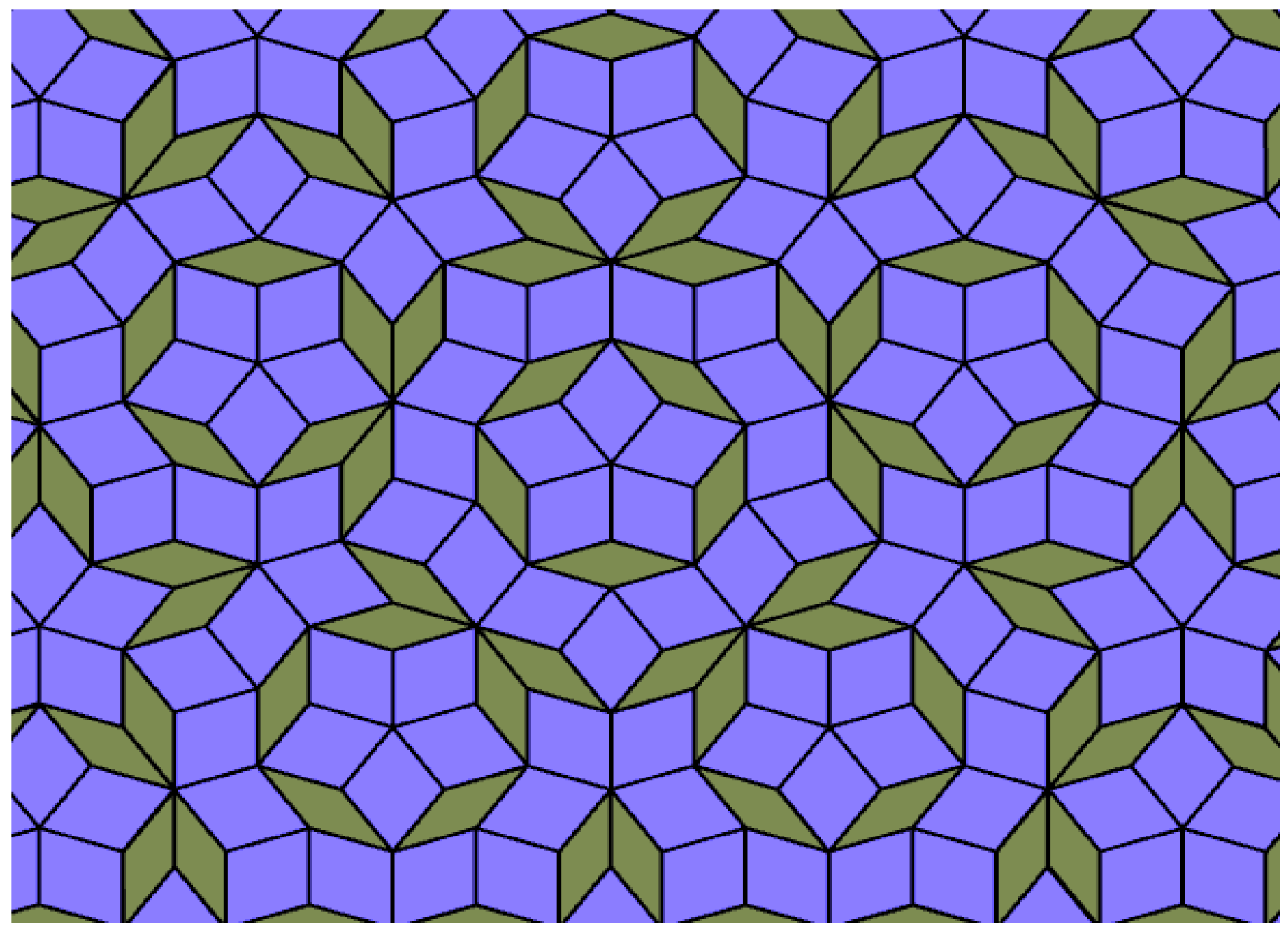

**Figure 9**. Excerpt of a Penrose tessellation with the fat and thin rhombuses.

On Figure 8 (1, far left), for example, five fat rhombuses are combined and five "horns" (outer vertices) are featured.

Each of the two rhombuses (thin or fat) and each of the five rhombus combinations (see Figure 8) is represented by a motif with a number of measures. The main motif A, represented by 1 in Figure 8, appears in the tessellation surrounded by 2, 3 and 4 in Figure 8. It is musically made up of a 5/8 measure motif to represent the five "horns" (the 2D view of the image), as well as three additional measures. This motif begins on the first voice (mandolin voice) and is imitated by the guitar at measure 4, creating a superimposition of the motif. Motif B is represented by Figure 2 (two measures in 5/8). Motif C (see 3 in Figure 8) consists of three 5/8 measures. The final motif D (see 4 in Figure 8) is composed of two 5/8 measures. For their first appearance, motifs A-B-C-D are consecutive for the mandolin voice and for the guitar voice (which begins at measure 4). The harp will present motif B at measures 5 and 6 (superimposed with motif B in the mandolin voice) and at measures 9 and 10.

We can model the entries and times of the instruments (mandolin, harp or guitar) in the first movement by a set of sequences of the form (Start time, Duration, Motif), where Start time and Duration are positive integers and Motif is either A, B, C or D.

The sequence for the mandolin will then be:

$S_{mandolin} = \{(0, 3, A), (3, 2, B), (5, 3, C), (8, 2, D), (10, 3, A), (13, 2, B), (15, 3, C)\ldots\}$.

For the guitar, which imitates the mandolin's motifs starting from measure 4, we have:

$S_{guitar} = \{(3, 3, A), (6, 2, B), (8, 3, C), (11, 2, D), (13, 3, A), (16, 2, B), (18, 3, C)\ldots\}$.

The harp also follows a sequence with overlaps:

$S_{harp} = \{(4, 2, B), (8, 2, B), (12, 2, B)\ldots\}$.

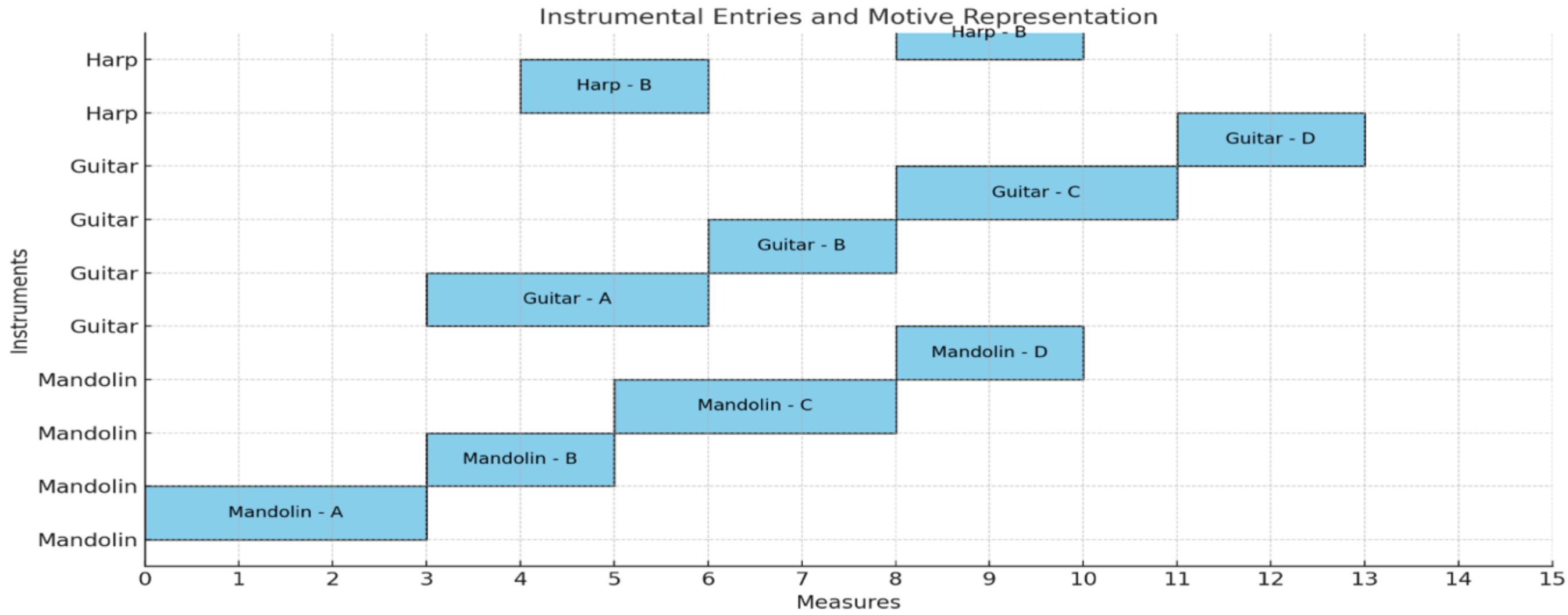

**Figure 10**. Instrumental entries on *Tesselles sonores* (first movement).

The guitar imitates the mandolin after three measures and plays the same motifs, in the same order. The cycle then repeats itself after ten measures; however, for aesthetic reasons, the composer has chosen not to stick rigorously to the repeating process and has taken some artistic licence with the full cycle after the ten first measures. Musically, as it has already been emphasised in the text and in Figure 10, the duration of motifs A and C is three measures, whereas motifs B and D last two measures. For each instrument, it will be convenient to represent the sequences of start times and associated durations as a data matrix, with $t_i$ the start times and $d_i$ the durations:

| **Instrument** | **Motif** | **Start time ($t_i$)** | **Duration ($d_i$)** |
|---|---|---|---|
| Mandolin | *A* | 0 | 3 |
| Mandolin | *B* | 3 | 2 |
| Mandolin | *C* | 5 | 3 |
| Mandolin | *D* | 8 | 2 |
| Guitar | *A* | 3 | 3 |
| Guitar | *B* | 6 | 2 |
| Guitar | *C* | 8 | 3 |
| Guitar | *D* | 11 | 2 |
| Harp | *B* | 4 | 2 |
| Harp | *B* | 8 | 2 |

The timing function that describes the overlap of each motive is defined by the following elementary counting function:

$$T(t) = \sum_i \mathbf{1}_{\{t_i \le t \le t_i + d_i\}}$$

where **1** is the indicator function, that equals 1 if the time $t$ is within the start time $t_i$ and the end time $t_i + d_i$ for that motive, and 0 otherwise.

For a given motif $M$, the overlap function between two instruments $I_1$ and $I_2$ can be defined as the union of the intersections of the time intervals $[t_j, t_j + d_j]$ when $I_1$ is playing $M$ with the time intervals $[t_k, t_k + d_k]$ when $I_2$ is playing $M$:

$$\text{Overlap}(I_1, I_2, M) = \bigcup_{j,k} \left( \left[ t_j, t_j + d_j \right] \cap \left[ t_k, t_k + d_k \right] \right).$$

Take for example the second motif ($M = B$) and $I_1$ = Mandolin, $I_2$ = Harp. Both instruments play motif $B$ during the interval [4, 5], then during the intervals [13, 14], [24, 25]…

The pattern thus creates a repetitive and aperiodic structure, just like a Penrose tiling, where the sequence does not repeat exactly, but follows a base structure that interrelates and transforms itself (similar to how the Penrose tiles fit together without repeating).

For the second and third movements, the composer selected the famous Nasrid Bird tessellation (see Figure 11), which can be seen in the Nasrid palaces of the Alhambra in Granada, Andalusia (Spain).

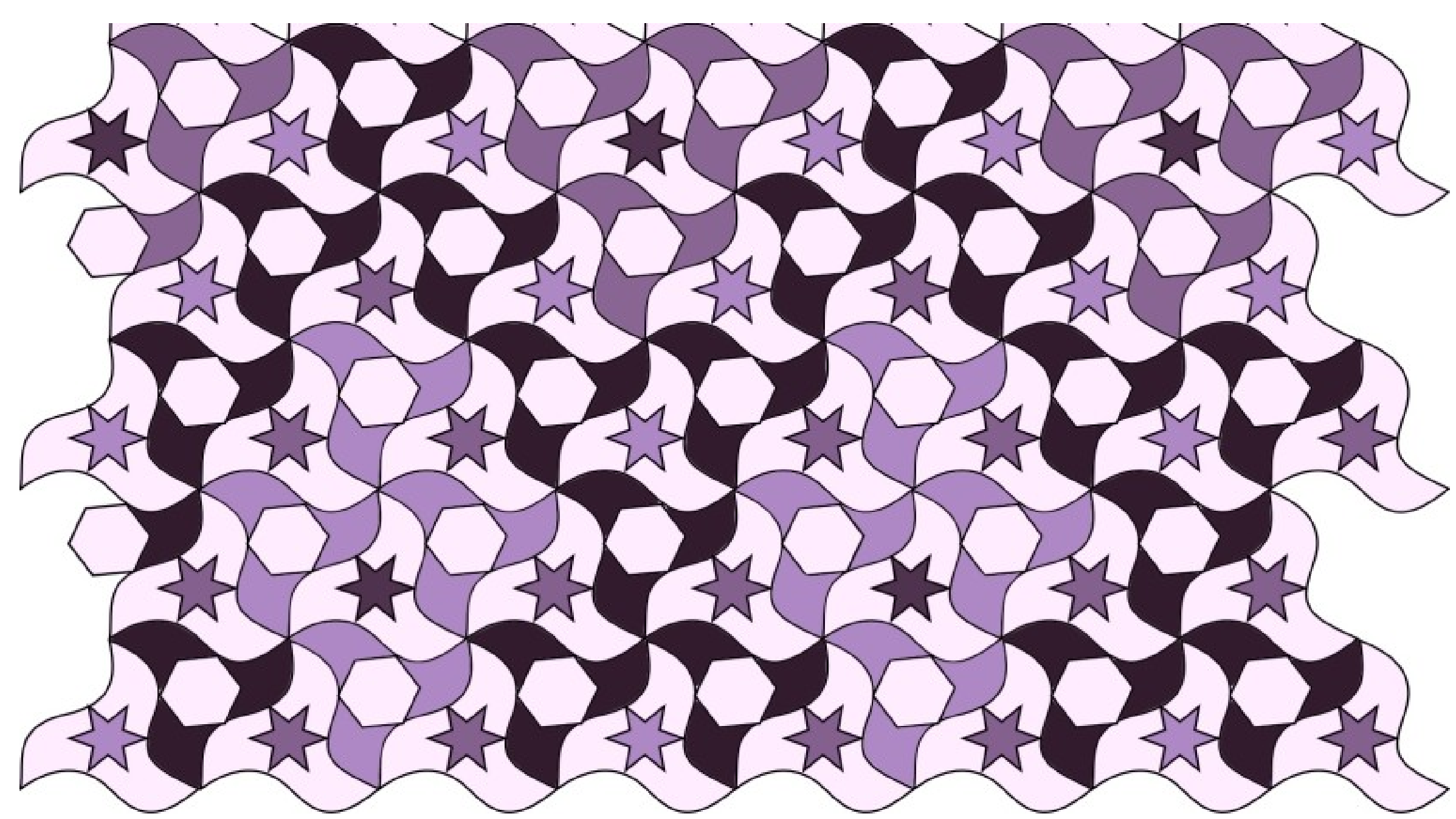

**Figure 11**. Representation of the Nasrid Bird tessellation.

Let us name each colour of the mosaic with a letter, allowing for a row matrix (or vector) encoding of the different colour combinations. The mosaic particularly features two distinct elements: the "stars" and the "birds" (see Figure 12).

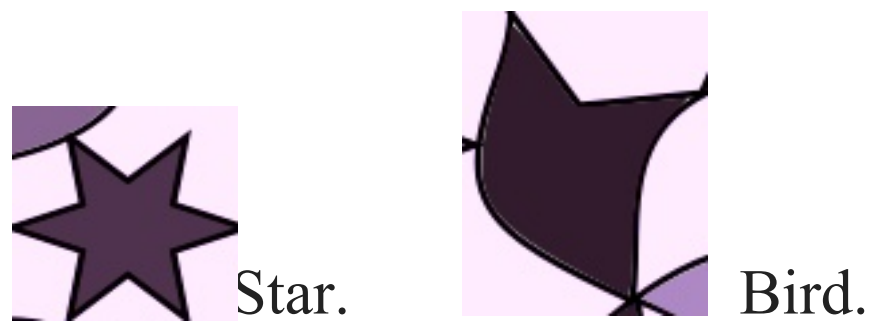

**Figure 12**. Representation of a star (left) and a bird (right) from the mosaic.

We will assign a shade and a letter to each of the tessellation's colours in order to distinguish them. There are three different colours of stars (see Figure 11): A (dark), B (medium), and C (light). Since these are isolated elements, the matrices for each colour will simply be chosen as follows: (1) for A, (2) for B, (1 2) for C. The motifs for A and B will consist of a single measure; the motifs for C will consist of two measures.

The birds are always presented in groups of three around a star (see Figure 11). The birds also appear in three colours (see Figure 11): A (dark, like some stars), B (medium, like some stars), and a new colour D (white). The matrices representing each group of three birds will be respectively: (1 1 1) for A, (2 1 1) for B, and (1 2 2) for D. Each group consists of three measures.

For the needs of the composition, the mosaic was divided into six parts. Each part (detailed below) consists of three groups of three birds of the same colour surrounding a coloured star. The mosaic presents six combinations of colours (see Figure 11):

Part I

D(1, 2, 2)

D (1, 2, 2) A (1) B (2,1,1)

Part II

A (1, 1, 1)

A (1, 1, 1) C (1, 2) A (1, 1, 1)

Part III

D (1, 2, 2)

D (1, 2, 2) A (1) D (1, 2, 2)

Part IV

Dx3 (1, 2, 2)

D (1, 2, 2) B (2) A (1, 1, 1)

Part V

B (2, 1, 1)

B (2, 1, 1)    C (1,2)    B (2, 1, 1)

Part VI

D (1, 2, 2)

A (1, 1, 1)    B (2)    A (1, 1, 1)

All the possible combinations of the birds and stars that appear in the mosaic have been described and coded; see also Figure 13. Musically speaking, the motifs corresponding to the six parts appear simultaneously, distributed among the three instruments.

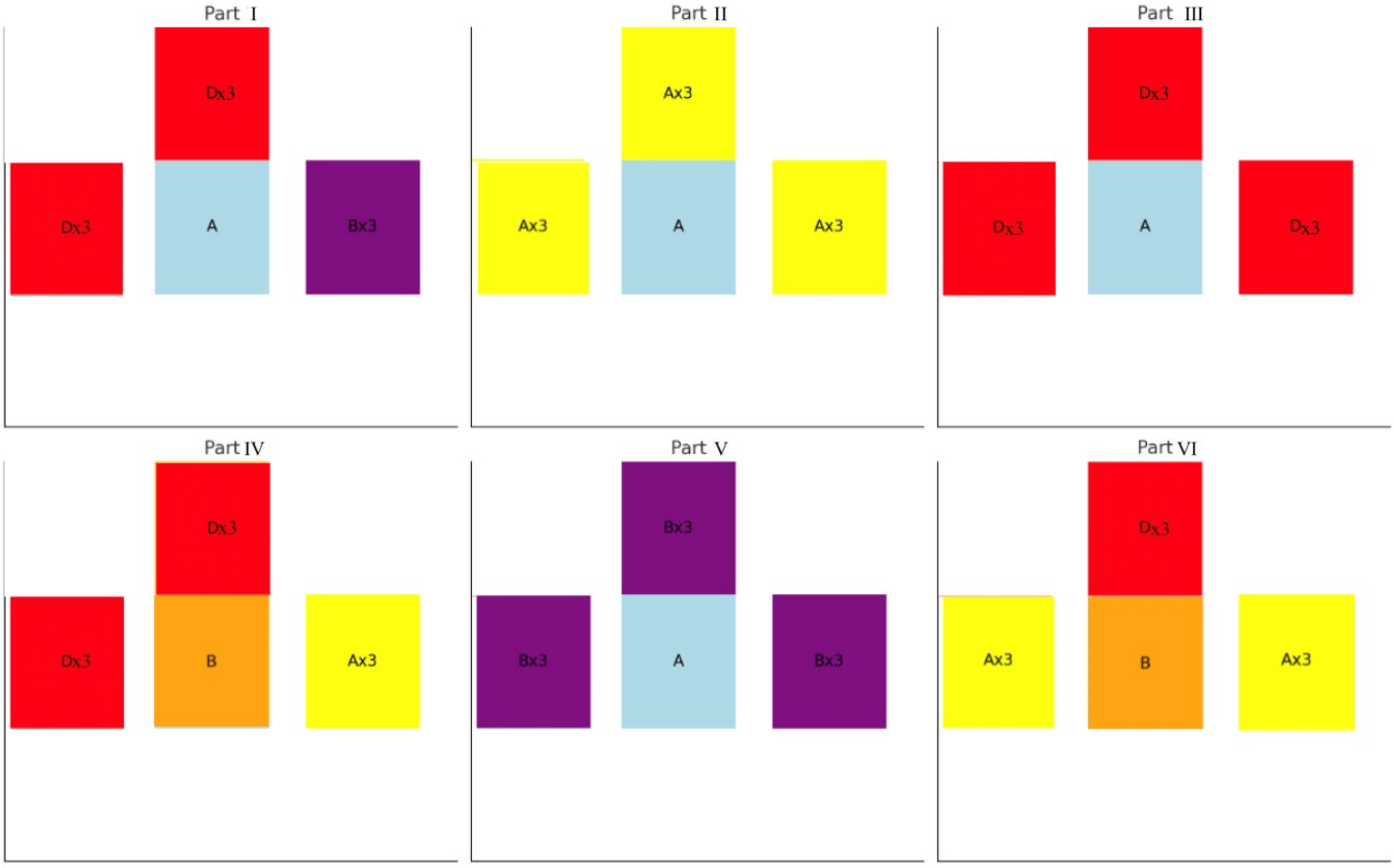

**Figure 13**. Colour distribution of the mosaic.

## 4. Parallel cases: mosaics and Steve Reich. Similarities with *Clapping Music*

There are other similarities between the canon entries of the isorhythmic motet, created through rhythmic shifting, and certain contemporary works. The minimalist piece *Clapping Music* by Steve Reich (1972) presents an example of using rhythmic shifting to create musical structures that lack exact repetition but are nonetheless regular. The piece was composed for two percussionists performing the same rhythmic motif. It relies on a principle of phase shifting, where one percussionist begins to play the motif earlier than the other, thus creating a phenomenon of progressive rhythmic shifting. This process generates a sound structure that seems to evolve without repeating, marked by aperiodicity, even though the total duration of the cycle remains constant.

The rhythmic dynamics of *Clapping Music* find a connection with recent research on aperiodic tilings, such as that presented by Fujita and Niizeki, who proposes a new variant of the famous Penrose rhombic tiling (see Figure 9). This tiling, although aperiodic, is governed by strict assembly rules for the two types of rhombuses, thus ensuring coverage of the plane without periodic repetition. Their work introduces a modified version of Penrose rhombic tiling, maintaining the fundamental principles of aperiodicity while adding elements that facilitate its geometric understanding. On the one hand, the tiling, famous for using only two distinct shapes, involves symmetries of order 5. (Fujita and Niizeki, 2025) introduce a cycle of ten units, significantly altering the pattern and relationships between the tiles. The idea of a cycle of ten units might have an interesting parallel with the concept of repetition and rhythmic shifting in *Clapping Music*. Although Reich works with shorter cycles (the complete cycle is twelve measures in total), the idea that the cycles "dilate" and overlap creates a sense of rhythmic aperiodicity similar to what occurs with Fujita and Niizeki's tiling, thus opening avenues for exploration for composers.

For performer A, let us assume that the base pattern $P_A$ has a duration $T$ and is divided into $n$ parts (generally, $n = 12$): the pattern can be expressed as a sequence of $n$ values corresponding to the durations of each clap. We thus have $P_A = [p_1, p_2 \dots p_n]$, where $p_i$ represents the time at which the $i^{th}$ clap occurs.

The second performer (B) begins by clapping the same pattern, but he shifts progressively. In other words, B's pattern $P_B$ starts with a shift $\Delta t > 0$ relative to $P_A$. After each cycle, the shift increases by a constant amount. This can be represented as $P_B(t) = [p_1 + \Delta t, p_2 + \Delta t \dots p_n + \Delta t]$, where $\Delta t$ gradually increases, typically at a constant fraction of the complete cycle $T$.

To model the gradual change in the shift, we can assume that the shift increases linearly with time $t$, so that at the $k^{th}$ cycle, the shift is proportional to $\Delta t_k = \frac{k}{n} T$, where $T$ is the duration of the cycle and $k$ is the cycle number. Thus, the sequence for the second performer becomes:

$$P_B(k) = \left[ p_1 + \frac{k}{n} T, p_2 + \frac{k}{n} T \dots p_n + \frac{k}{n} T \right]$$

where $k$ is the number of cycles elapsed, and $(k/n)T$ is the shift, which increases with $k$.

## 5. Conclusion

David Smith's discovery opens up numerous artistic possibilities regarding the use of musical elements that can allow for the tiling of the plane in unexpected ways, not only in relation to the geometric transformation used in *Le Chapeau à douze cornes*, which transitions from symmetry to asymmetry: aperiodic tessellations from a single tile open up new ways of exploration for composers. Moreover, the Penrose tiling and those of Fujita and Niizeki −with some amplifications− in the works *Tesselles sonores* and *Clapping Music* show a similarity in the use of Vuza's canon. However, a few particularities are worth noting: Movement I of *Tesselles sonores* (with Penrose tilings) illustrates this idea of the canon; the concept behind Movements II and III blends combinatorics, permutations, and the translation of a set of elements. In the case of the isorhythmic motet example, two elements act in canon (talea and color) until the plane is filled, whereas in *Clapping Music*, the focus is on translation movements (as seen in Fujita and Niizeki's tilings) −only one voice is shifted.

As we have shown in this article, the geometric theory of tilings can be applied in connection with musical criteria such as pitch, duration, or sound material. The arrangement of sound material (following the visual aspects of tilings) offers a different perspective on the abstract conception of composition. At the same time −as we have observed− new modelling of tilings (both periodic and aperiodic related to those of Penrose) can be connected to historical compositional techniques (as in the case of the isorhythmic motet) and other more recent techniques (such as the phase shifting that Steve Reich employs in some of his pieces, like *Clapping Music*).

## Acknowledgements

We thank David Smith (the discoverer of the "Hat") for all his valuable information and comments on his explorations, as well as Craig Kaplan for all the discussions on this subject. Thanks to Anne Montaron and the team of the program Création Mondiale (France Musique-Radio France) for making the recording and popularization of the development of the piece *Le Chapeau à douze cornes* possible. The first authors would also like to thank the ensemble L'Instant Donné for the creation of this piece as well as the Trio Polycordes for their interpretation of *Tesselles sonores*. Their performances of the works motivated and inspired us to write this article.

## Disclosure statement

In this article, the artwork of the first author, Maria Luisa Acuña Fuentes, is referred to under her artist name, Marisa Acuña. No potential conflict of interest was reported by the authors. Figures 9 (excerpt from a Penrose tiling) and 11 (excerpt from a tessellation at the Alhambra) have been slightly adapted from unattributed public sources. All other figures have been entirely designed and crafted by the authors.

## Références